\documentclass[12pt]{article}
\RequirePackage{amsfonts,amssymb,amsbsy,amsmath}
\usepackage{epsfig,longtable}

\DeclareMathOperator{\diag}{diag}

\def\no{\noindent}

\newtheorem{defn}{Definition}[section]

\newtheorem{exa}{\sc \bf Example}[section]
\newtheorem{rem}{\sc Remark}[section]

\def \y {{\bf y}}

\def \RR {{\mathbb{R}}}

\def \no {\noindent}

\def \pmatrix{ \left( \begin{array} }
\def \endpmatrix{ \end{array} \right) }
\begin{document}

\date{}

\title{Fifty Years of Stiffness\thanks{Work developed
within the project ``Numerical methods and software for
differential equations''.}}

\author{Luigi Brugnano\thanks{Dipartimento di
Matematica, Universit\`a di Firenze, Viale Morgagni 67/A, 50134
Firenze (Italy).\qquad E-mail:\,{\tt luigi.brugnano@unifi.it}}
\and Francesca Mazzia\thanks{Dipartimento di Matematica,
Universit\`a di Bari, Via Orabona  4,  70125 Bari (Italy).\qquad
\mbox{E-mail:}\,{\tt mazzia@dm.uniba.it}} \and
 Donato Trigiante\thanks{Dipartimento di Energetica, Universit\`a di Firenze,
Via Lombroso 6/17, 50134 Firenze (Italy).\qquad E-mail:\,{\tt
trigiant@unifi.it}} }

\maketitle

\begin{abstract}
The notion of {\em stiffness}, which originated in several
applications of a different nature, has dominated the activities
related to the numerical treatment of differential problems for
the last fifty years. Contrary to what usually happens in
Mathematics, its definition has been, for a long time, not
formally precise (actually, there are too many of them). Again,
the needs of applications, especially those arising in the
construction of robust and general purpose codes, require nowadays
a formally precise definition. In this paper, we review the
evolution of such a notion and we also provide a precise
definition which encompasses all the previous ones.

\bigskip \noindent{\bf Keywords:} stiffness, ODE problems,
discrete problems, initial value problems, boundary value
problems, boundary value methods.
\end{abstract}

\section{Introduction}

\begin{flushright}
  \begin{minipage}[t]{3in}
     \textsf{\em Frustra fit per plura quod potest per pauciora.}
\begin{flushright}
    \textsf{Razor of W.\,of Ockham, {\em doctor invincibilis}.}
\end{flushright}
  \end{minipage}
\end{flushright}

The struggle generated by the duality short times--long times is
at the heart of  human culture in almost all its aspects. Here are
just a few examples to fix the idea:

\begin{itemize}
  \item in historiography:  Braudel's distinction among the geographic,
social and individual times;\footnote{\,Moreover, his concept of
\emph{structure}, i.e. events which are able to accelerate the
normal flow of time, is also interesting from our point of view,
because  it somehow  recalls the mathematical concept of large
variation in small intervals of time (see later).}

\item in the social sphere: Societies are organized according to three kinds of
laws, i.e., codes   (regulating short term relations), constitutions
(regulating medium terms   relations), and ethical laws (long term
rules) often not explicitly stated but religiously accepted;

  \item in the economy sphere: the laws of this part of human activities
  are partially unknown at the moment. Some models (e.g., the Goodwin model
\cite{good}), permits us to say, by taking into  account only a
few variables, that the main  evolution is periodic in time (and
then predictable), although we are experiencing an excess of
periodicity (chaotic behavior). Nevertheless, some
  experts claim (see, e.g., \cite{Ga}) that the problems in the predictability of
  the economy  are mainly due to a sort of gap in passing information
  from a generation to the next ones, i.e. to the conflict between
  short time and long time behaviors.\footnote{\,Even Finance makes the
  distinction between short time and long time traders.}
\end{itemize}

Considering the importance of this concept, it would have been
surprising if the duality ``short times--long times''  did not
appear somewhere in Mathematics. As a matter of fact, this
struggle not only appears in our field but it  also has a name:
\emph{stiffness}.

Apart from a few early papers \cite{CN,Hir}, there is a general
agreement in placing the date of the introduction of such problems
in Mathematics to around 1960 \cite{Da6}. They were the
necessities of the applications to draw the attention of the
mathematical community towards such problems, as the name itself
testifies: ``\emph{they have been termed stiff since they
correspond to tight coupling between the driver and the driven
components in servo-mechanism}" (\cite{Da1} quoting from
\cite{Hir}).

Both the number and the type of applications proposing difficult
differential problems has increased exponentially in the last
fifty years. In the early times, the problems proposed by
applications were essentially initial value problems and,
consequently, the definition of stiffness was clear enough and
shared among the few experts, as the following three examples
evidently show:

\begin{description}
  \item[D1]:
\emph{Systems containing very fast components as well as very slow
components} \\ (Dahlquist \cite{Da1}).

  \item[D2]: \emph{They represent coupled physical systems having
  components varying with very different times scales: that is they
  are systems having some components varying much more rapidly
  than the others} (Liniger \cite{Li}, translated from French).

  \item[D3]: \emph{A stiff system is one for which $\lambda_{max}$ is
  enormous so that either the stability or the error bound or both
  can only be assured by unreasonable restrictions on $h$\dots
  Enormous means enormous relative to the scale which here is
  $\bar{t}$ (the integration interval)\dots} (Miranker \cite{Mi}).
\end{description}

The above definitions are rather informal, certainly very far from
the precise definitions we are accustomed to in Mathematics, but,
at least, they agree  on a crucial point: the relation among
stiffness and the appearance  of different time-scales in the
solutions (see also \cite{Hind}).

 Later on, the necessity to encompass new classes of difficult problems, such as
Boundary Value Problems, Oscillating Problems, etc., has led
either to weaken the definition or, more often, to define some
consequence of the phenomenon instead of defining the phenomenon
itself.  In Lambert's book \cite{Lam} five propositions about
stiffness, each of them capturing some important aspects of it,
are given. As matter of fact, it has been also stated that no
universally accepted definition of stiffness exists \cite{schoST}.

There are, in the literature, other definitions based on other
numerical difficulties, such as, for example, large  Lipschitz
constants or logarithmic norms \cite{Sod},  or non-normality of
matrices \cite{TH}. Often is not even clear if stiffness refers to
particular solutions (see, e.g. \cite{Hu}) or to problems as a
whole.

Sometimes one has the feeling that stiffness is becoming so broad
to be nearly synonymous of difficult.

At the moment, even if the old intuitive definition relating
stiffness to multiscale problems survives in most of the authors,
the most successful definition seems to be the one
  based on  particular  effects of the phenomenon rather than
 on the phenomenon itself, such as, for example, the following
 almost equivalent items:

\begin{description}
\item[D4]: \emph{Stiff equations are equations where certain implicit
methods \dots perform better, usually tremendous better, than
explicit ones} \cite{Hir}.

  \item[D5]: \emph{Stiff equations are problems for which explicit
  methods don't work} \cite{Hairer}.

  \item[D6]: \emph{If a numerical method with a finite region of absolute
  stability, applied to a system with any initial condition, is
  forced to use in a certain interval of integration a step length
  which is excessively small in relation to the smoothness of the
  exact solution in that interval, then the system is said to be
  stiff in that interval} \cite{Lam}.
\end{description}

As usually happens, describing a phenomenon by means of  its
effects may not be enough to fully characterize the phenomenon
itself. For example, saying that fire is what produces ash, would
oblige firemen to wait for the end of a fire to see if the ash has
been produced. In the same way, in order to recognize stiffness
according to the previous definitions, it would be necessary to
apply first one\footnote{\,It is not clear if one is enough: in
principle the definition may require to apply  all of them.}
explicit method and see if it works or not.
Some authors, probably discouraged by the above defeats in giving a rigorous
definition, have also affirmed that a rigorous
mathematical definition of stiffness is not possible
\cite{schoGH}.

It is clear that this situation is unacceptable for at least two
reasons:

\begin{itemize}
  \item it is against the tradition of Mathematics, where objects
  under study have to be \emph{precisely} defined;
  \item it is necessary to have the possibility to recognize \emph{operatively} this
    class of problems, in order to increase the efficiency of the
    numerical codes to be used in applications.
\end{itemize}

Concerning the first item, our opinion is that, in order to gain
in \textrm{precision}, it would be necessary to revise the concept
of \emph{stability} used in Numerical Analysis, which is somehow
different from the homonym concept used in all the other fields of
Mathematics, where stable are equilibrium points, equilibrium
sets, reference solutions, etc., but not equations or
problems\footnote{\,Only in particular circumstances, for example
in the linear case, it is sometimes allowed the language abuse:
the nonlinear case may contain simultaneously  stable and unstable
solutions.} (see also \cite{Da6} and \cite{LeVeq}).

Concerning the second item, \emph{operatively} is intended in the
sense that the definition must be stated in terms of
\emph{numerically observable} quantities such as, for example,
norms of vectors or matrices. It was believed that, seen from the
applicative point of view, a formal definition of stiffness would
not be strictly necessary: \emph{Complete formality here is of
little value to the scientist or engineer with a real problem to
solve} \cite{Hind}.

Nowadays, after the great advance in the quality of numerical
codes,\footnote{\,A great deal of this improvement is due to the
author of the previous sentence.} the usefulness of a formal
definition is strongly recognised, also from the point of view of
applications: \emph{One of the major difficulties associated with
the study of stiff differential systems is that a good
mathematical definition of the concept of stiffness does not
exist} \cite{Ca}.

 In this paper, starting from ideas already partially exposed
elsewhere \cite{bt2,bt1,imt}, we will try to unravel the
question of the definition of stiffness and show that a precise
and operative definition of it, which encompasses all the
known facets, is possible.

In order to be as clear as possible, we shall start with the
simpler case of initial value for a single linear equation and
gradually we shall consider more general cases and, eventually, we
shall synthesize the results.

\section{The asymptotic stability case}

For initial value problems for ODEs, the concept of stability
concerns the behavior of a generic solution $y(t)$, in the
neighborhood of a reference solution $\bar{y}(t)$, when the
initial value is perturbed. When the problem is linear and
homogeneous, the difference, $e(t) = y(t)-\bar{y}(t)$, satisfies
the same equation as $\bar{y}(t)$. For nonlinear problems, one
resorts to the linearized problem, described by the variational
equation, which, essentially, provides valuable information only
when $\bar{y}(t)$ is asymptotically stable. Such a variational
equation can be used to generalize to nonlinear problems the
arguments below which, for sake of simplicity, concerns only the
linear case.

Originally, stiffness was almost always associated with initial
value problems having asymptotically stable equilibrium points
(dissipative problems) (see, e.g., Dahlquist \cite{Da2}). We then
start from this case, which is a very special one. Its
peculiarities arise from the following two facts:\footnote{\,We
omit, for simplicity, the other fact which could affect new
definitions, i.e., the fact that the solutions of the linear
equation can be integrated over any large interval because of the
equivalence, in this case, between asymptotic and exponential
stability.}

\begin{itemize}
  \item it is the most common in applications;
  \item there exists a powerful and fundamental theorem, usually called
\emph{Stability in the first approximation
  Theorem} or \emph{Poincar\'{e}-Liapunov Theorem}, along with its corollary due
to Perron\footnote{\,It is interesting to observe that
  the same theorem is known as the {\em Ostrowsky's Theorem}, in the theory of
iterative methods.}, which allows us
  to reduce the study of stability of critical points, of a very large
  class of nonlinearities, to the study of the stability of the corresponding
  linearized problems (see, e.g., \cite{Co,LL,RM,Yo}).
\end{itemize}

The former  fact explains the pressure of applications for the
treatment of such problems even before the computer age. The
latter one provides, although not always explicitly recognized,
the mathematical solid bases for the profitable and extensive use,
in Numerical Analysis, of the linear test equation to study the
fixed-$h$ stability of numerical methods.

We shall consider explicitly the case where the linearized problem
is autonomous, although the following definitions will take into
account the more general case.

Our starting case will then be that of an initial value problem
having an asymptotically stable reference solution, whose
representative is, in the scalar case,

\begin{eqnarray}\label{single}
 y'&=&\lambda y,  \qquad
 t\in[0,T], \qquad \mbox{Re}\lambda <0,\\
 y(0)&=&\eta,\nonumber
\end{eqnarray}

\no where the reference solution (an equilibrium point, in this
case) has been placed at the origin. From what is said above, it
turns out that it is not by chance that it coincides with the
famous test equation.

\begin{rem}
It is worth observing that the above test equation is not less
general than $y'=\lambda y+g(t)$, which very often appears in the
definitions of stiffness: the only difference is the reference
solution, which becomes $\bar{y}(t)=\int_0^t e^{\lambda
(t-s)}g(s)ds$, but not the topology of solutions around it. This
can be easily seen by introducing the new variable
$z(t)=y(t)-\bar{y}(t)$ which satisfies exactly  equation
(\ref{single}) and then, trivially,  must  share the same
stiffness. Once the solution $z(t)$ of the homogeneous equation
has been obtained, the solution $y(t)$ is obtained by adding to it
$\bar{y}(t)$ which, in principle, could be  obtained by  means of
a quadrature formula. This allows us to conclude that if any
stiffness is in the problem, this must reside in the homogeneous
part of it, i.e., in problem (\ref{single}).
\end{rem}

\begin{rem}
We call attention to the interval of integration $[0,T]$, which
depends on our need for information about the solution, even if
the latter exists for all values of $t$. This interval must be
considered as datum of the problem. This has been sometimes
overlooked, thus creating some confusion.
\end{rem}

Having fixed problem (\ref{single}), we now look for a
mathematical tool which allows us to state formally the intuitive
concept, shared by  almost all the definitions of stiffness: i.e.,
we look for one or two parameters which tells us if in $[0,T]$ the
solution varies rapidly or not. This can be done easily by
introducing the following two measures for the solution of problem
(\ref{single}):

\begin{equation}\label{k1k2}
  \kappa_c=\frac{1}{|\eta|}\max_{t\in[0,T]} |y(t)|, \qquad
  \gamma_c=\frac{1}{|\eta|}\frac{1}{T}\int_0^T|y(t)|dt,
\end{equation}

\no which, in the present case, assume the values:
$$
\kappa_c=1, \qquad \gamma_c=\frac{1}{|Re \lambda |T}\,(1-e^{Re
\lambda  T})\approx \frac{1}{|Re\lambda|T}=\frac{T^*}{T},
$$
\no where $T^*=|Re\lambda|^{-1}$ is the transient time. The two measures
$\kappa_c$, $\gamma_c$ are called \emph{conditioning
parameters} because they measure the sensitivity of the solution subject to a
perturbation of the initial conditions in the infinity and in the $l_1$ norm.

Sometimes, it would be preferable to use a lower value of
$\gamma_c$, i.e.,

\begin{equation}\label{alternat}
  \gamma_c=\frac{1}{|\lambda |T}.
\end{equation}

\no This amounts to consider also  the oscillating part of
the solution (see also Remark~\ref{frequency} below).

\begin{figure}
\centerline{\includegraphics[width=8cm,height=8cm]{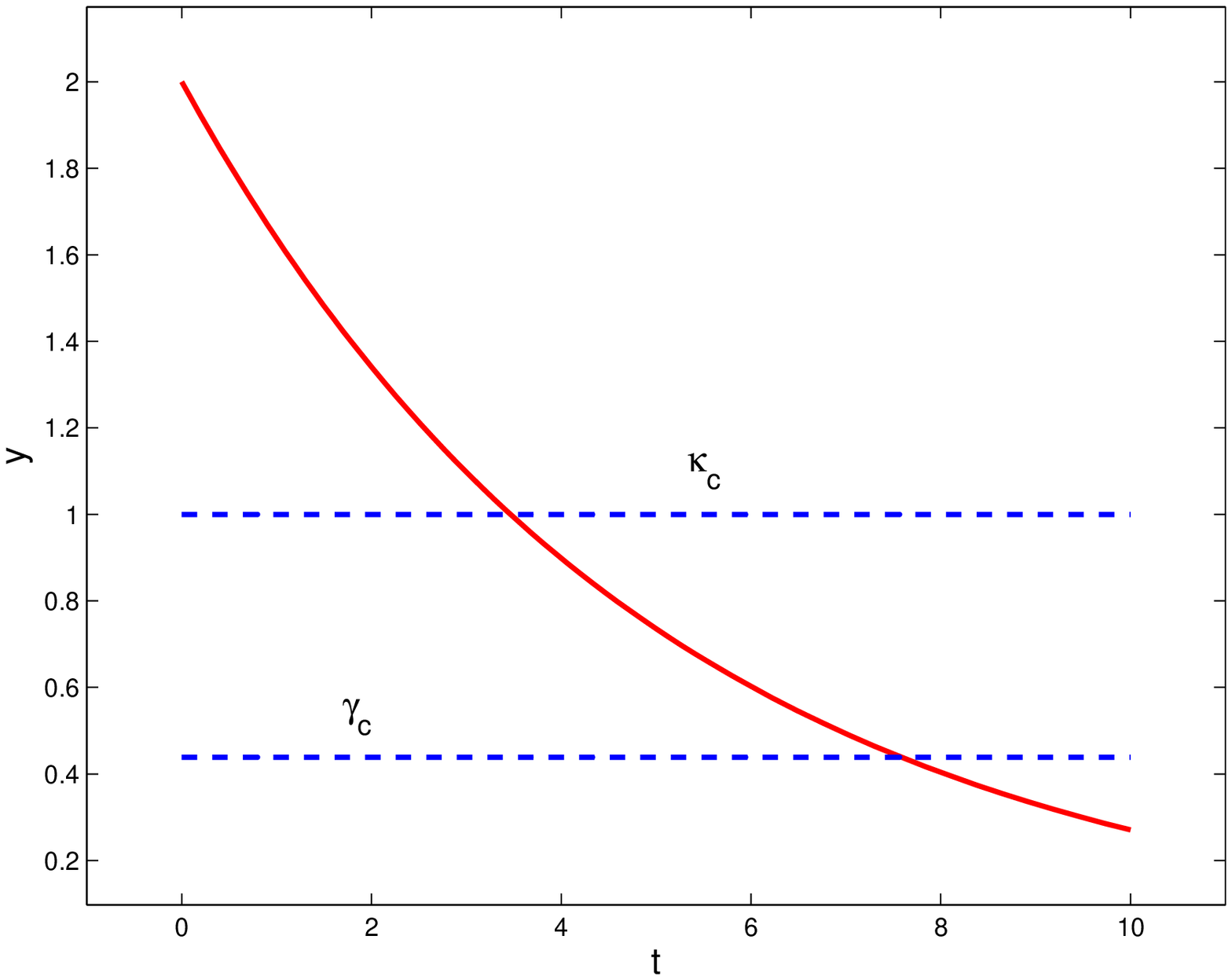}
\includegraphics[width=8cm,height=8cm]{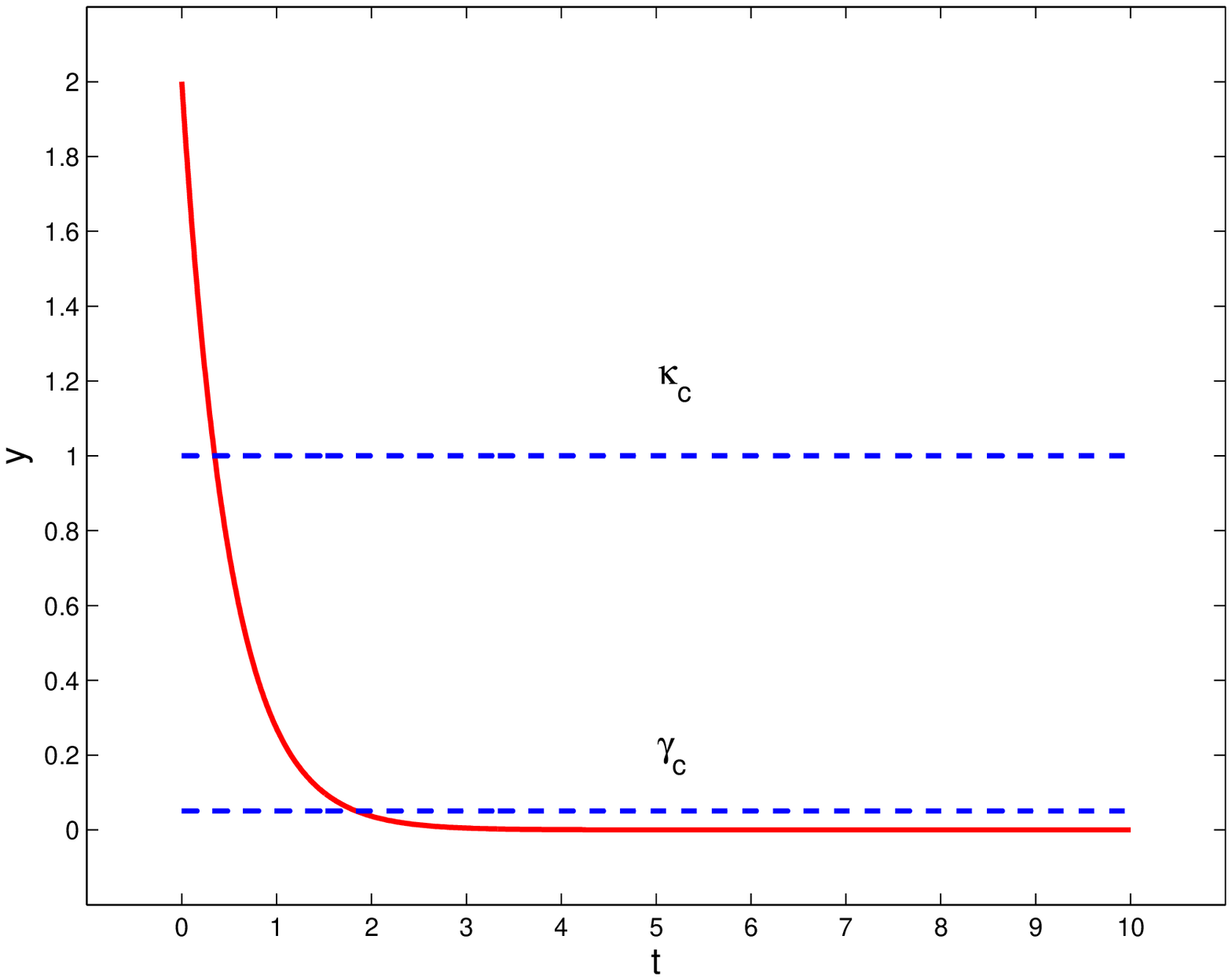}}
\caption{\protect\label{kgam} Solutions and values of $k_c$ and $\gamma_c$ in
the cases $\lambda=-.2$ (left plot) and $\lambda=-2$ (right plot).}
\end{figure}

By looking at Figure\,\ref{kgam}, one realizes at once that a
rapid variation of the solution in $[0,T]$ occurs when
$k_c\gg\gamma_c.$ It follows then that the parameter

\begin{equation}\label{sigma}
  \sigma_c=\frac{k_c}{\gamma_c}\equiv \frac{T}{T^*},
\end{equation}

\no which is the ratio between the two characteristic times of the
problem, is more significant. Consequently, the definition of
stiffness follows now trivially:

\begin{defn}\label{stiffsingle}
The initial value problem (\ref{single}) is {\em stiff} if
$\sigma_c\gg1$.
\end{defn}

The parameter $\sigma_c$ is called \emph{stiffness ratio}.

\begin{rem}
The width of the integration interval $T$ plays a fundamental role
in the definition. This is an important point: some authors, in
fact, believe that  stiffness should concern equations; some
others believe that stiffness should concern problems, i.e.,
equations and data. We believe that both statements are partially
correct: stiffness concerns equations,  integration time, and a
set of initial data (not a specific one of them). Since this point
is more important in the non scalar case, it will be discussed in
more detail later.
\end{rem}

\begin{rem}\label{frequency}
When $\gamma_c$ is defined according to (\ref{alternat}), the definition of
stiffness continues to be also meaningful in the case $Re \lambda=0$, i.e., when
the critical point is only marginally stable. In fact, let
$$\lambda=i\omega\equiv i\frac{2\pi}{T^*}.$$ Then,
$$\sigma_c=2\pi\frac{T}{T^*},$$ and the definition
encompasses also the case of \emph{oscillating stiffness}
introduced by some authors (e.g., \cite{Mi}). Once again the
stiffness is the ratio of two times. If information about the
solution on the smaller time scale is needed, an adequately small
stepsize should be used. It is worth noting that high oscillating
systems (with respect to $T$) fall in the class of problems for
which explicit methods do not work, and then are stiff according
to definitions D4--D6.

When $\lambda=0$, then $k_c=\gamma_c=\sigma_c=1.$

In the case $Re\lambda >0$ (i.e., the case of an unstable critical
point), both parameters $k_c$ and $\gamma_c$ grow exponentially
with time. This implies that small variations in the initial
conditions will imply exponentially large variations in the
solutions, both pointwise and on average: i.e., the problem is
\emph{ill conditioned}.
\end{rem}

Of course, the case $Re\lambda = 0$ considered above cannot be
considered as representative of more difficult nonlinear equations,
since linearization is in general not allowed in such a case.

The linearization is not the only way to study nonlinear
differential (or difference) equations. The so called {\em
Liapunov second method} can be used as well (see, e.g.,
\cite{Hahn,LL,Yo}). It has been used, in connection with stiffness
in \cite{Bu,Da2,Da3,Da4,Da5,Da6}, although not always explicitly
named.\footnote{\,Often, it appears under the name of one-sided
Lipschitz condition.}    Anyway, no matter how the  asymptotic
stability of a reference solution is detected, the parameters
(\ref{k1k2}) and Definition~\ref{stiffsingle} continue to be
valid. Later on, the problem of effectively estimating such
parameters will also be discussed.

\subsection{The discrete case}
Before passing to the non scalar case, let us now consider the
discrete case, where some interesting additional considerations
can be made. Here, almost all we have said for the continuous case
can be repeated. The first approximation theorem can be stated
almost in the same terms as in the continuous case (see e.g.
\cite{LakTri}).

Let the interval $[0,T]$ be partitioned into $N$ subintervals of
length $h_n>0$, thus defining the mesh points: $t_n = \sum_{j=1}^n
h_j,$ $n=0,1,\ldots,N.$

The linearized autonomous problem is now:

\begin{equation}\label{singledis}
  y_{n+1}=\mu_n y_n, \qquad n=0,\dots, N-1, \qquad y_0=\eta,
\end{equation}

\no where the $\{\mu_n\}$ are complex parameters. The conditioning
parameters  for (\ref{singledis}), along with the stiffness ratio,
are defined as:

\begin{equation}\label{par}
  \kappa_d=\frac{1}{|\eta|}\max_{i=0,\dots,N}|y_i|,\qquad \gamma_d=
  \frac{1}{|\eta|}\frac{1}{T}\sum_{i=1}^{N}h_i \max(|y_i|,|y_{i-1}|),\qquad
  \sigma_d=\frac{k_d}{ \gamma_d}.
\end{equation}

This permits us to define the notion of {\em well representation}
of a continuous problem by means of a discrete one.

\begin{defn}\label{wr}
The problem (\ref{single}) is {\em well represented} by
(\ref{singledis}) if
\begin{eqnarray}\label{wr1}
k_c&\approx &k_d,\\ \label{wr2} \gamma_c&\approx& \gamma_d.
\end{eqnarray}
\end{defn}

In the case of a constant mesh-size $h$, $\mu_n\equiv \mu$ and it
easily follows that the condition (\ref{wr1})  requires $|\mu|<1$.
It is not difficult to recognize the usual $A$-stability
conditions for one-step methods (see Table~\ref{tb1}).
Furthermore, it is easily recognized that the request that
condition (\ref{wr1}) holds uniformly with respect to $h\lambda\in
\mathbb{C}^-$  implies that the numerical method producing
(\ref{singledis}) must be implicit.

\begin{table}
  \centering
  \begin{tabular}{|c|c|c|}\hline
    method & $\mu$ & condition \\ \hline
    Explicit Euler& $1+h\lambda$ & $|1+h\lambda|<1$ \\[1mm]
    Implicit Euler & $\frac{1}{1-h\lambda}$ &
$\left|\frac{1}{1-h\lambda}\right|<1$\\[2mm]
    Trapezoidal Rule  & $\frac{1+h\lambda/2}{1-h\lambda/2}$ &
    $\left|\frac{1+h\lambda/2}{1-h\lambda/2}\right|<1$ \\ \hline
  \end{tabular}
  \caption{\protect\label{tb1} Condition (\ref{wr1}) for some popular
methods.}
\end{table}

What does condition (\ref{wr2}) require more? Of course, it
measures how faithfully the integral $ \int_0^T|y(t)|dt$ is
approximated by the quadrature formula $\sum_{i=1}^Nh_i
\max(|y_i|,|y_{i-1}|)$, thus giving a sort of global information
about the behavior of the method producing the approximations
$\{y_i\}$. One of the most efficient global strategies for
changing the stepsize is based on monitoring this parameter
\cite{bt3,bt1,CM,CMT,MT1,MT2}. In addition to this, when finite
precision arithmetic is used, then an interesting property of the
parameter $\gamma_d$ occurs \cite{imt}: if it is smaller than a
suitably small threshold, this suggests that we are doing useless
computations, since the machine precision has already been
reached.

\subsection{The non scalar case}

In this case, the linearized problem to be considered is

\begin{equation}\label{vector}
  y'=Ay, \qquad t\in[0,T],\qquad y(0)=\eta,
\end{equation}

\no with $A\in\RR^{m\times m}$ and having all its eigenvalues with
negative real part. It is clear from what was said in the scalar
case that, denoting by $\Phi(t)=e^{At}$ the fundamental matrix of
the above equation, the straightforward generalization of the
definition of the conditioning parameters (\ref{k1k2}) would lead
to:

\begin{equation}\label{nonsingle}
  \kappa_c=\max_{t\in[0,T]}\|\Phi(t)\|, \qquad
  \gamma_c=\frac{1}{T}\int_0^T\|\Phi(t)\|dt, \qquad
\sigma_c=\frac{\kappa_c}{\gamma_c}.
\end{equation}

\no Indeed, these straight definitions \emph{work most of the
time}, as is confirmed by the following example, although, as we
shall explain soon, not always.

\begin{figure}
\centerline{\includegraphics[width=12cm,height=8cm]{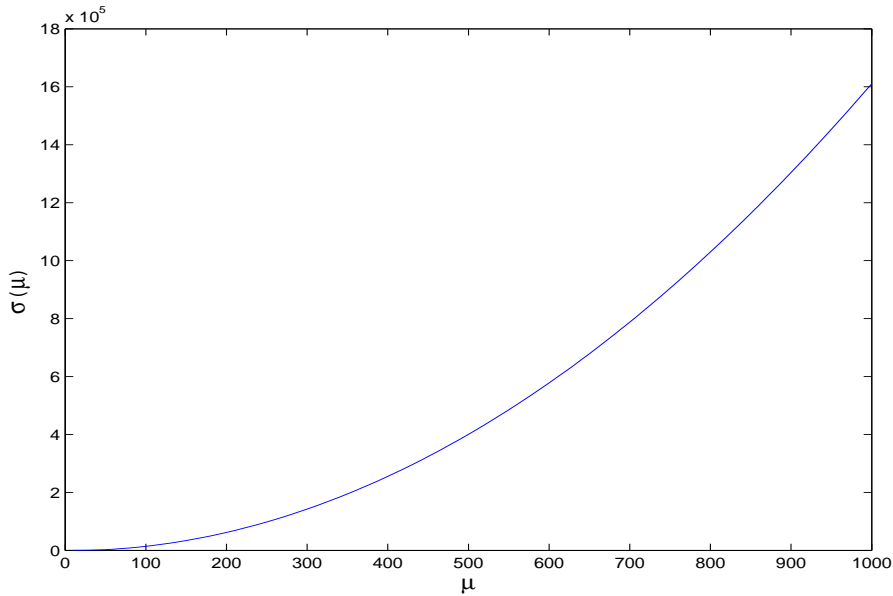}}
\caption{\protect\label{vdp} Estimated stiffness ratio of Van der Pol's
problem (\ref{vdpol}).}
\end{figure}

\begin{exa} Let us consider the well-known Van der Pol's problem,
\begin{eqnarray}\nonumber
y_1' &=& y_2, \\ \label{vdpol}
y_2' &=& -y_1+\mu\, y_2(1-y_1^2), \qquad t\in[0,2\mu], \\
y(0) &= & (2,~0)^T,\nonumber
\end{eqnarray}

\no whose solution approaches a limit cycle of period $T\approx 2
\mu$. It is also very well-known that, the larger the parameter
$\mu$, the more difficult the problem is. In Figure~\ref{vdp} we
plot the parameter $\sigma_c(\mu)$ (as defined in
(\ref{nonsingle})) for $\mu$ ranging from 0 to $10^3$. Clearly,
stiffness increases with $\mu$.
\end{exa}

Even though (\ref{nonsingle}) works for this problem, this is not
true in general. The problem is that the definition of stiffness
as the ratio of two quantities may require a lower bound for the
denominator. While the definition of $\kappa_c$ remains unchanged,
the definition of $\gamma_c$ is more entangled. Actually, we need
two different estimates of such a parameter:

\begin{itemize}
  \item an upper bound,  to be used for estimating the
  conditioning  of the problem in $l_1$ norm;
  \item a lower bound, to be used in defining $\sigma_c$ and, then,
  the stiffness.
\end{itemize}

In the  definition given in \cite{bt2,bt1}, this distinction was
not made, even though the definition was (qualitatively) completed
by adding
\begin{equation}\label{quot}
\mbox{``\emph{for at least one of the modes}''.}
\end{equation}
\no We shall be more precise in a moment. In the meanwhile, it is interesting to
note that the clarification contained in (\ref{quot}) is already in one
of the two definitions given by Miranker \cite{Mi}:

\emph{A system of differential equations is said to be stiff on
the interval $(0,\bar{t})$ if there exists a solution of that
system a component of which has a variation on that interval which
is large compared to $\frac{1}{\bar{t}}$},

\no where it should be stressed that the definition considers
equations and not problems: this   implies that the existence of
largely variable components may appear for at least one choice of
the initial conditions, not necessary for a specific one.

Later on, the definition was modified so as to translate into
formulas the above quoted sentence (\ref{quot}). The following
definitions were then given (see, e.g., \cite{imt}):

\begin{equation}
\label{parcont}
\begin{array}{ll}
\kappa_c(T,\eta)=\displaystyle \frac{1}{\|\eta\|} \max_{0 \le t
\le T}\|y(t)\|, &\qquad \kappa_c(T)=\displaystyle \max_{\eta}
{\kappa_c(T,\eta)},
\\~\\
\gamma_c(T,\eta)=
 \displaystyle \frac{1}{T \|\eta\|}\int_{0}^{T} \|y(t) \|dt,  & \qquad
\gamma_c(T)=\displaystyle \max_{\eta} {\gamma_c(T,\eta)}.
\end{array}
\end{equation}

\no and

\begin{equation} \label{parcont1}
\sigma_c(T) =\max_{\eta}
\frac{\kappa_c(T,\eta)}{\gamma_c(T,\eta)}.
\end{equation}

The only major change regards the definition of $\sigma_c.$ Let us
be more clear on this point with an example, since it leads to a
controversial question  in the literature: i.e., the dependence of
stiffness from the initial condition.  Let
$A=\diag(\lambda_1,\lambda_2,\ldots,\lambda_m)$ with $\lambda_i<0$
and $|\lambda_1|>|\lambda_2|>\ldots>|\lambda_m|$. The solution of
problem (\ref{vector}) is $y(t)=e^{At}\eta$.

If  $\sigma_c$  is defined according to (\ref{nonsingle}), it
turns out that $\|e^{At}\|=e^{\lambda_m t}$ and, then,
$\gamma_c(T)\approx \frac{1}{T|\lambda_m|}.$ If, however, we take
$\eta=(1,0,\ldots,0)^T$, then $y(t)=e^{\lambda_1t}$ and
$\gamma_c(T)$ becomes $\gamma_c(T)\approx \frac{1}{T|\lambda_1|}$.
Of course, by changing the initial point, one may activate each
one of the \emph{modes}, i.e. the functions $e^{\lambda_it}$ on
the diagonal of the matrix $e^{At}$, leaving silent the others.
This is the reason for specifying, in the older definition, the
quoted sentence (\ref{quot}). The new definition (\ref{parcont1}),
which essentially poses as the denominator of the ratio $\sigma_c$
the smallest value among the possible values of
$\gamma_c(T,\eta)$, is more compact and complies with the needs of
people working on the construction of codes, who like more
operative definitions. For the previous diagonal example, we have
that $k_c$ continues to be equal to 1, while
$\gamma_c(T)=\frac{1}{T|\lambda_1|}$.

Having got the new definition (\ref{parcont1}) of $\sigma_c(T)$,
the definition of stiffness continues to be given by
Definition~\ref{stiffsingle} given in the scalar case, i.e., the
problem (\ref{vector}) is {\em stiff} if\, $\sigma_c (T)\gg1$.

How does this definition reconcile with the most used definition
of stiffness for the linear case, which considers the ``smallest"
eigenvalue  $\lambda_m$ as well? The answer is already in
Miranker's definition D3. In fact, usually the integration
interval is chosen large enough  to provide complete information
on the behavior of the solution. In this case, until the slowest
mode has decayed enough, i.e. $T=1/|\lambda_m|$, which implies

\begin{equation}\label{stiffmat}
  \sigma_c \left (T=\frac{1}{|\lambda_m|} \right )
  =\left|\frac{\lambda_1}{\lambda_m}\right|,
\end{equation}

\no which, when much larger than 1, coincides with the most
common definition of stiffness in the linear case. However, let us
insist on saying that if the interval of integration is much
smaller than $1/|\lambda_m|$, \emph{the problem may be not stiff} even
if $\left|\frac{\lambda_1}{\lambda_m}\right|\gg1$.

The controversy about the dependence of the definition of
stiffness on the initial data is better understood by considering
the following equation given in \cite[pp.\,217--218]{Lam}:

$$
\frac{d}{dt}\pmatrix{c}y_1\\y_2\endpmatrix=\pmatrix{cc}-2&1\\-1.999&0.999\endpmatrix
\pmatrix{c}y_1\\y_2\endpmatrix+\pmatrix{c}2\sin t\\0.999(\sin
t-\cos t )\endpmatrix,
$$

\no whose general solution is

$$
\pmatrix{c}y_1\\y_2\endpmatrix=c_1e^{-t}\pmatrix{c}1\\1
\endpmatrix+c_2e^{-0.001t}\pmatrix{c}1\\1.999\endpmatrix+\pmatrix{c}\sin t\\\cos
t\endpmatrix.
$$

\no The initial condition $y(0)=(2,\,3)^T$ requires $c_2=0$ and,
then, the slowest mode is not activated: the solution rapidly
reaches the reference solution. If this information was known
beforehand, one could, in principle,  choose  the interval of
integration $T$ much smaller than $\frac{1}{0.001}$. This,
however, does  not take into account the fact that the computer
uses finite precision arithmetic, which may not represent exactly
the initial condition $\eta$. To be more precise, let us point out
that the slowest mode is not activated only if the initial
condition is on the line $y_2(0)-y_1(0)-1=0$. Any irrational value
of $y_1(0)$ will not be well represented on the computer. This is
enough to activate the silent mode. Of course, if one is sure that
the long term contribution to the solution obtained on the
computer is due to this kind of error, a small value of $T$ can
always be used. But it is rare that this information is known in
advance. For this reason, we consider the problem to be stiff,
since we believe  that  the definition of stiffness cannot
distinguish, for example,  between rational and irrational values
of the initial conditions. Put differently, initial conditions are
like a fuse that may activate stiffness.

We conclude this section by providing a few examples, which show
that Definition~\ref{stiffsingle}, when $\sigma_c$ is defined
according to (\ref{parcont1}), is able to adequately describe the
stiffness of nonlinear and/or non autonomous problems as well.

\begin{figure}[t]
\centerline{\includegraphics[width=12cm,height=8cm]{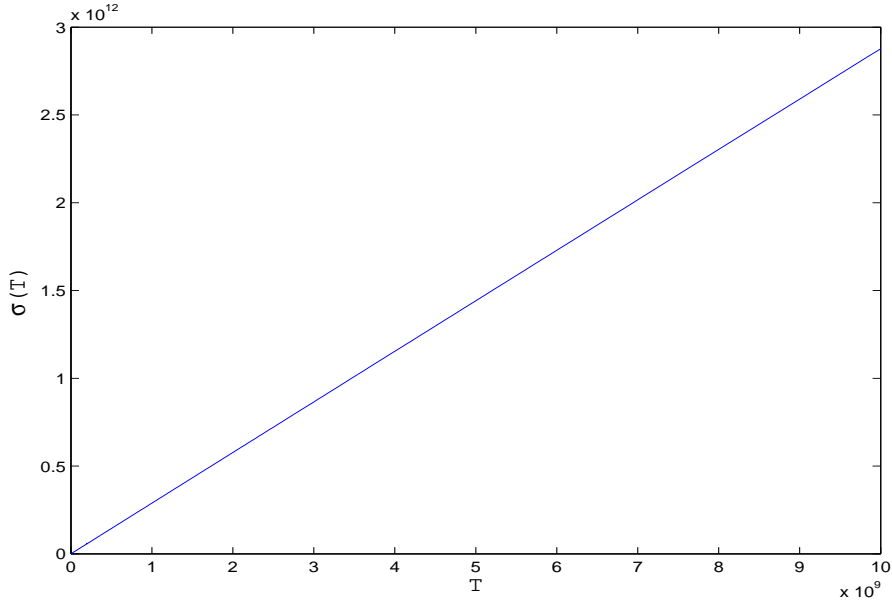}}
\caption{\protect\label{rober} Estimated stiffness ratio of Robertson's
problem (\ref{rob}).}
\end{figure}

\begin{exa}
Let us consider the well-known Robertson's problem:
\begin{eqnarray}\nonumber
  y_1' &=& -.04y_1 +10^4y_2y_3, \\ \label{rob}
  y_2' &=&  .04y_1 -10^4y_2y_3 -3\cdot 10^7y_2^2,\qquad t\in[0,T],\\
\nonumber
  y_3' &=&                      3\cdot 10^7y_2^2,\\
\nonumber
   y(0) &=& (1,~0,~0)^T.
\end{eqnarray}

\no Its stiffness ratio with respect to the length $T$ of the
integration interval, obtained through the linearized problem and
considering a perturbation of the initial condition of the form
$(0,\,\varepsilon,\,-\varepsilon)^T$, is plotted in
Figure~\ref{rober}. As it is well-known, the figure confirms that
for this problem stiffness increases with $T$.
\end{exa}

\begin{figure}[t]
\centerline{\includegraphics[width=12cm,height=8cm]{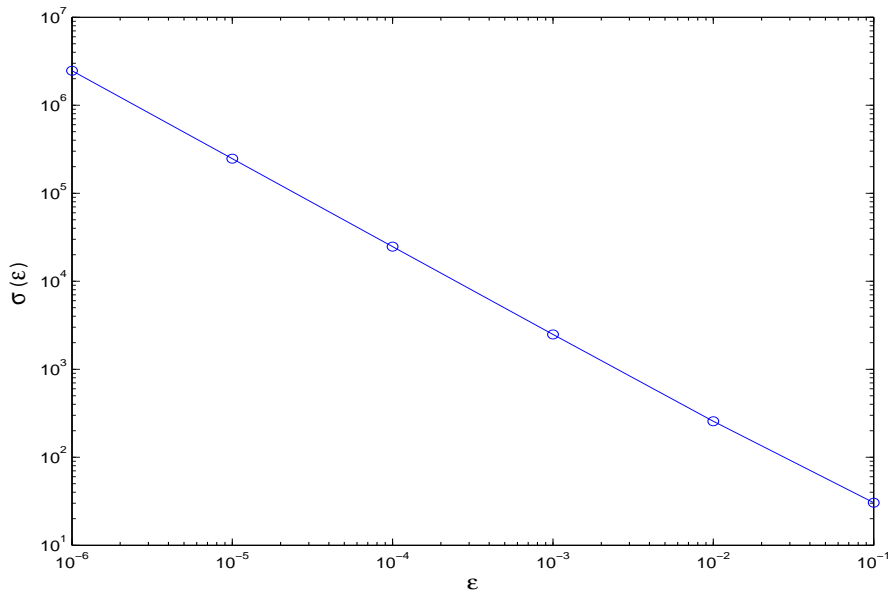}}
\caption{\protect\label{kreiss_hairer} Estimated stiffness ratio of the
Kreiss problem (\ref{kreh})--(\ref{kreh2}).}
\end{figure}

\begin{exa}
Let us consider the so-called Kreiss problem
\cite[p.\,542]{Hairer}, a linear and non autonomous problem:

\begin{equation}\label{kreh}
y' = A(t) y, \qquad t\in[0,4\pi], \qquad y(0) \mbox{\quad fixed},
\end{equation}

\no where
\begin{equation}\label{kreh1}
A(t) =  Q^{T}(t) \Lambda_\varepsilon Q(t),
\end{equation}

\no and
\begin{equation}\label{kreh2}
 Q(t) =
\pmatrix{cc} \cos t & \sin t \\ -\sin t & \cos t \endpmatrix,
\qquad \Lambda_\varepsilon = \pmatrix{cc} -1\\ &
-\varepsilon^{-1}\endpmatrix.
\end{equation}

\no Its stiffness ratio with respect to the small positive
parameter $\varepsilon$, obtained by considering a perturbation of
the initial condition of the form $(-\varepsilon,\,1)^T$, is
plotted in Figure~\ref{kreiss_hairer}. As one expects, the figure
confirms that the stiffness of the problem behaves as
$\varepsilon^{-1}$, as $\varepsilon$ tends to 0.
\end{exa}

\begin{figure}[ht]
\centerline{\includegraphics[width=12cm,height=8cm]{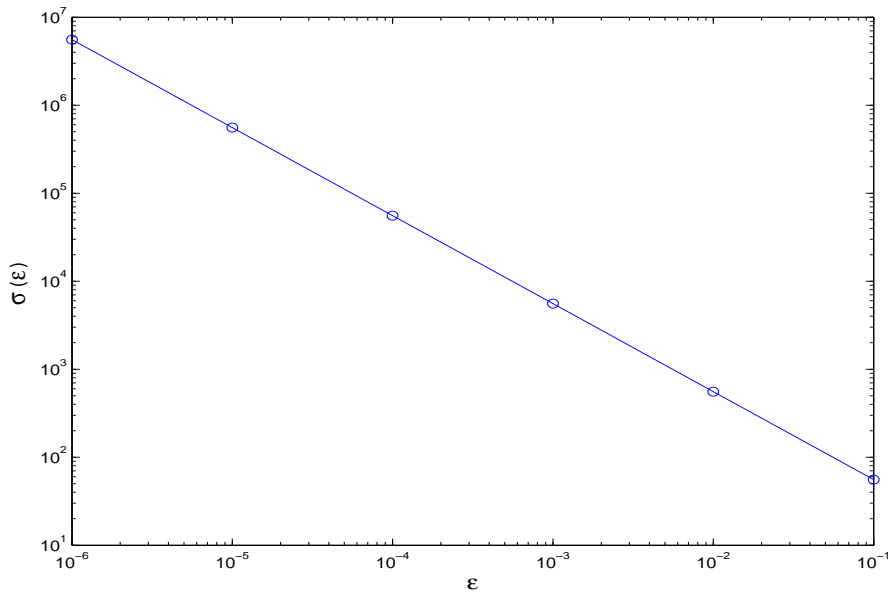}}
\caption{\protect\label{kreiss} Estimated stiffness ratio of the modified
Kreiss problem (\ref{kre})--(\ref{kre2}).}
\end{figure}

\begin{exa}
Let us consider the following linear and non autonomous problem, a
modification of problem (\ref{kreh}), that we call ``modified
Kreiss problem'':\,\footnote{\,This problem has been suggested by
J.I.\,Montijano.}

\begin{equation}\label{kre}
y' = A(t) y, \qquad t\in[0,4\pi], \qquad y(0) \mbox{\quad fixed},
\end{equation}

\no where
\begin{equation}\label{kre1}
A(t) = Q_\varepsilon^{-1}(t) P^{-1} \Lambda_\varepsilon P
Q_\varepsilon(t),
\end{equation}

\no and
\begin{equation}\label{kre2}
P = \pmatrix{rr} -1& 0\\\ 1&1\endpmatrix, \qquad Q_\varepsilon(t)
= \pmatrix{cc} 1& \varepsilon \\ e^{\sin t} & e^{\sin
t}\endpmatrix, \qquad \Lambda_\varepsilon = \pmatrix{cc} -1\\ &
-\varepsilon^{-1}\endpmatrix.
\end{equation}

\no Its stiffness ratio with respect to the small positive
parameter $\varepsilon$, obtained by considering a perturbation of
the initial condition of the form $(-\varepsilon,\,1)^T$, is shown
in Figure~\ref{kreiss}. Also in this case the stiffness of the the
problem behaves as $\varepsilon^{-1}$, as $\varepsilon$ tends to
0.
\end{exa}

\begin{rem}
It is worth mentioning that, in the examples considered above, we
numerically found that
$$
\max_\eta\, \frac{\kappa_c(T,\eta)}{\gamma_c(T,\eta)}
$$
is obtained by considering an initial condition $\eta$ in the
direction of the eigenvector of the Jacobian matrix (computed for
$t \approx t_0$) associated to the dominant eigenvalue. We note
that, for an autonomous linear problem, if $A$ is diagonalizable,
this choice activates the {\em mode} associated with $\lambda_1$,
i.e., the eigenvalue of maximum modulus of $A$.

\end{rem}

\subsection{The non scalar discrete case}

As for the scalar case, what we said for the continuous problems
can be repeated, {\em mutatis mutandis}, for the discrete ones.
For brevity, we shall skip here the details for this case, also
because they can be deduced from those described in the more
general case discussed in the next section.

\section{Boundary Value Problems (BVPs)}\label{BVPsec}

The literature about BVPs is far less abundant than that about
IVPs, both in the continuous and in the discrete case. While there
are countless books on the latter subject presenting it from many
points of view (e.g., stability of motion, dynamical systems,
bifurcation theory, etc.), there are many less books about the
former. More importantly, the subject is usually presented as a by
product of the theory of IVPs. This is not necessarily the best
way to look at the question, even though many important  results
can be obtained this way. However, it may sometimes be more useful
to look at the subject the other way around. Actually, the
question is that IVPs are naturally  a subclass of BVPs.   Let us
informally clarify this point without many technical details which
can be found, for example, in \cite{bt1}.

IVPs transmit  the initial information ``from left to right''.
Well conditioned IVPs are those for which the initial value, along
with the possible initial errors, decay moving from left to right.
FVPs (Final Value problems) are those transmitting information
``from right to left'' and, of course, well conditioning should
hold when the time, or the corresponding independent variable,
varies towards $-\infty$. More precisely, considering the scalar
test equation (\ref{single}), the asymptotically stability for
IVPs and FVPs requires $Re\lambda <0$ and $Re \lambda>0$,
respectively. BVPs transmit information both ways. Consequently,
they cannot be scalar problems but vectorial of dimension at least
two. We need then  to refer to the test equation (\ref{vector}).
It can be affirmed that a well conditioned linear BVP needs to
have eigenvalues  with both negative and positive real parts
(\emph{dichotomy}, see, e.g., \cite{amr,bt1}). More precisely: the
number of eigenvalues with negative real part has to match the
amount of information transmitted ``from left to right'', and the
number of eigenvalues with positive real part has to match the
amount of information traveling ``from right to left''. For
brevity, we shall call the above statement \emph{continuous
matching rule}. Of course, if there are no final conditions, then
the problem becomes an IVP and, as we have seen, in order to be
well conditioned, it must have all the eigenvalues with negative
real part. In other words, the generalization of the case of
asymptotically stable IVPs is the class of well conditioned BVPs
\emph{because both satisfy the continuous matching rule}. This is
exactly what we shall assume hereafter.

Similar considerations apply to the discrete problems, where the
role of the imaginary axis is played by the unit circumference in
the complex plane. It is not surprising that a numerical method
will {\em well represent} a continuous autonomous linear BVP if
the corresponding matrix has as many eigenvalues inside the unit
circle as the number of initial conditions and as many eigenvalues
outside the unit circle as the number of final conditions
(\emph{discrete matching rule}).

\begin{rem} The idea that IVPs are a subset of BVPs is  at the root of
the class of methods called {\em Boundary Value Methods (BVMs)}
which permits us, thanks to the discrete matching rule, to define
high order and perfectly $A$-stable methods (i.e., methods having
the imaginary axis separating the stable and unstable domains),
which \emph{overcome the Dahlquist's barriers}, and are able to
solve both IVPs and BVPs (see, e.g., \cite{bt1}).\end{rem}

\begin{rem} From this point of view, the popular {\em shooting
method}, consisting of transforming a BVP into an IVP and then
applying a good method \emph{designed for IVPs}, does not appear
to be such a good idea. As matter of fact, even a very well
conditioned linear BVP, i.e. one which satisfies the continuous
matching rule, will be transformed in a badly conditioned IVP,
since the matrix of the continuous IVP shall, of course, contain
eigenvalues with positive real part. This will prevent the
discrete matching rule to hold.\end{rem}

\subsection{Stiffness for BVPs}

Coming back to our main question, stiffness for BVPs is now
defined by generalizing the idea already discussed in the previous
sections.

As in the previous cases, we shall refer to linear problems, but
the definitions will also be applicable to nonlinear problems as
well. Moreover, according to what is stated above, we shall only
consider the case where the problems are well conditioned (for the
case of ill conditioned problems, the arguments are slightly more
entangled, see e.g. \cite{CM}). Then, let us consider the linear
and non autonomous BVP:

\begin{equation}\label{contprob}
y' = A(t) y, \qquad  t\in[0,T],\qquad  B_0 y(0)+B_1 y(T) = \eta,
\end{equation}

\noindent where $y(t),\eta \in \RR^m$ and $A(t),B_0,
B_1\in\RR^{m\times m}$. The solution of the problem
(\ref{contprob}) is
$$
y(t) = \Phi(t)Q^{-1} \eta,
$$

\no where $\Phi(t)$ is the fundamental matrix of the problem such
that $\Phi(0)=I$, and $Q= B_a + B_b \Phi(T)$, which has to be
nonsingular, in order for (\ref{contprob}) to be
solvable.\footnote{\,Observe that, in the case of IVPs, $B_0=I$
and $B_1=O$, so that $Q=I$.}

As in the continuous IVP case, the conditioning parameters are
defined (see (\ref{parcont})) as:

\begin{eqnarray}
\nonumber \kappa_c(T,\eta)=\frac{1}{\|\eta\|} \max_{0 \le t \le
T}\|y(t)\|, &\quad& \kappa_c(T)=\displaystyle \max_{\eta}
\kappa_c(T,\eta),
\\ \label{parcond} \\ \nonumber
\gamma_c(T,\eta) =
 \frac{1}{T\|\eta\|} \int_0^T \|y(t) \|dt,  &
 \quad& \gamma_c(T)= \max_\eta
 \gamma_c(T,\eta).
\end{eqnarray}

\no Consequently, the stiffness ratio is defined as (see
(\ref{parcont1})):
$$
\sigma_c(T) =\max_{\eta}
\frac{\kappa_c(T,\eta)}{\gamma_c(T,\eta)},
$$

\no and the problem is stiff if $\sigma_c(T)\gg1$. Moreover, upper
bounds of $\kappa_c(T)$ and $\gamma_c(T)$ are respectively given
by:
\begin{equation}\label{upperbound}
\kappa_c(T) \le \max_{0 \le t \le T} \| \Phi(t)Q^{-1}\|, \qquad
\gamma_c(T) \le \frac{1}T \int_0^T \| \Phi(t)Q^{-1}\| dt.
\end{equation}

Thus, the previous definitions naturally extend to BVPs the
results stated for IVPs. In a similar way, when considering the
discrete approximation of (\ref{contprob}), for the sake of
brevity provided by a suitable one-step method over a partition
$\pi$ of the interval $[0,T]$, with subintervals  of length $h_i$,
$i=1,\dots,N$, the discrete problem will be given by
\begin{equation}\label{disbvp}
  y_{n+1}=R_n y_n, \qquad n=0,\dots, N-1, \qquad B_0
  y_0+B_1y_N=\eta,
\end{equation}

\no whose solution is given by $$y_n = \left(\prod_{i=0}^{n-1}
R_i\right) Q_N^{-1}\eta, \qquad Q_N = B_0 +B_1 \prod_{i=0}^{N-1}
R_i.$$

The corresponding discrete conditioning parameters are then
defined by:
\begin{eqnarray}
\nonumber \kappa_d(\pi,\eta)=\frac{1}{\|\eta\|} \max_{0 \le n \le
N}\|y_n\|, &\quad& \kappa_d(\pi)=\displaystyle \max_{\eta}
\kappa_d(\pi,\eta),
\\ \label{dparcond} \\ \nonumber
\gamma_d(\pi,\eta) =
 \frac{1}{T\|\eta\|} \sum_{i=1}^N h_i\max(\|y_i\|,\|y_{i-1}\|),  &
 \quad& \gamma_d(\pi)= \max_\eta
 \gamma_d(\pi,\eta),
\end{eqnarray}
and
$$
\sigma_d(\pi) = \max_{\eta}
\frac{\kappa_d(\pi,\eta)}{\gamma_d(\pi,\eta)}.
$$

According to Definition~\ref{wr}, we say that the discrete
problem\footnote{\,It is both defined by the used method and by
the considered mesh.} (\ref{disbvp}) {\em well represents} the
continuous problem (\ref{contprob}) if
\begin{equation}\label{wrc}
\kappa_d(\pi) \approx \kappa_c(T), \qquad \gamma_d(\pi) \approx
\gamma_c(T).\end{equation}

\begin{rem} It is worth mentioning that innovative mesh-selection
strategies for the efficient numerical solution of stiff BVPs have
been defined by requiring the match (\ref{wrc}) (see, e.g.,
\cite{bt3, bt1, CM, CMT, imt}).\end{rem}

\begin{figure}[t]
\centerline{\includegraphics[width=12cm,height=8cm]{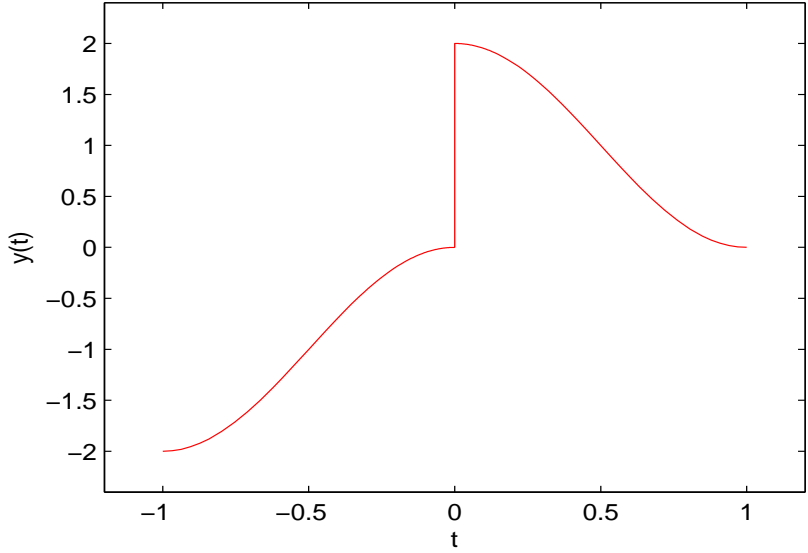}}
\caption{\protect\label{figbvp1f}  Problem (\ref{prob6}),
$\varepsilon=10^{-8}$.}
\end{figure}

\subsection{Singular Perturbation Problems}

The numerical solution of singular  perturbation  problems can be
very difficult because they can  have solutions   with very narrow
regions of rapid variation characterized by
 boundary layers, shocks, and interior layers. Usually, the equations
 depend on a small parameter, say $\varepsilon$,  and the problems become
 more difficult as $\varepsilon$ tends to 0.
 It is not always clear, however,
 how the width of the region of rapid variation is related to the parameter $\varepsilon$.
 By computing the stiffness ratio $\sigma_c(T)$, we observe that singularly
 perturbed problems are stiff problems. Moreover, as the following
 examples show,   the parameter $\sigma_c(T)$ provides us also with information
 about the width of the region of rapid variation.

 The examples are formulated as second order equations: of course, they have
 to be transformed into corresponding first order systems, in order to apply the
 results of the previous statements.

\begin{figure}[t]
\centerline{\includegraphics[width=12cm,height=8cm]{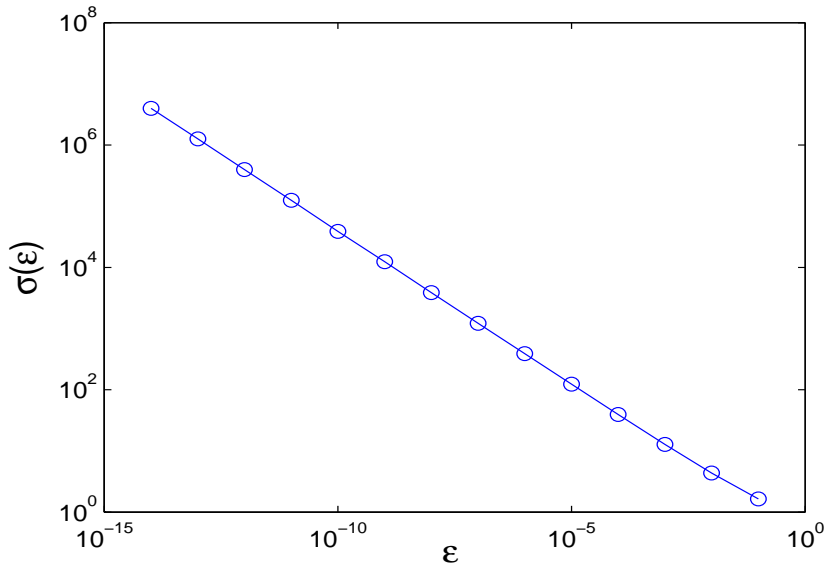}}
\caption{\protect\label{figbvp1} Estimated stiffness ratio of problem
(\ref{prob6}).}
\end{figure}

\begin{exa}
Let us consider the linear singularly perturbed problem:
\begin{equation}\label{prob6}
\begin{array}{l}
\varepsilon y'' + t y' = -\varepsilon \pi^2 \cos(\pi t) - \pi t
\sin(\pi t), \qquad y(-1) = -2, \quad y(1)=0,
\end{array}
\end{equation}
whose solution has, for $0 < \varepsilon \ll 1$, a turning point
at $t=0$ (see Figure \ref{figbvp1f}). The exact solution is
 $
  y(t)= \cos(\pi t)+ \exp((t-1)/\sqrt{\varepsilon})+\exp(-(t+1)/\sqrt{\varepsilon}).
$

In Figure \ref{figbvp1} we plot an estimate of the stiffness ratio
obtained by considering two different perturbations of the
boundary conditions of the form $(1, 0)^T$ and $(0,1)^T$. The
parameter $\varepsilon$ varies from $10^{-1}$ to  $10 ^{-14}$. We
see that the (estimated) stiffness parameter grows like
$\sqrt{\varepsilon^{-1}}$.
\end{exa}

\begin{figure}[t]

\centerline{\includegraphics[width=12cm,height=8cm]{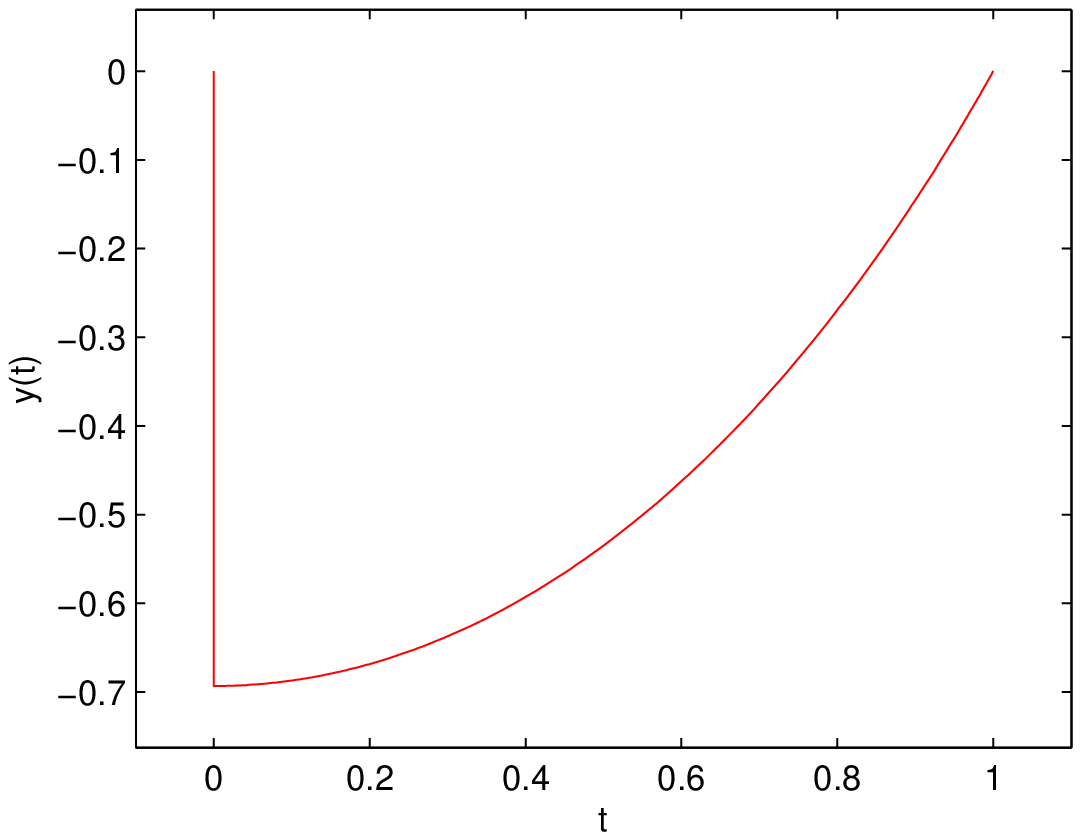}}
\caption{\protect\label{figbvp2f}  Problem (\ref{prob19}),
$\varepsilon=10^{-6}$.}
\end{figure}

\begin{figure}[t]
\centerline{\includegraphics[width=12cm,height=8cm]{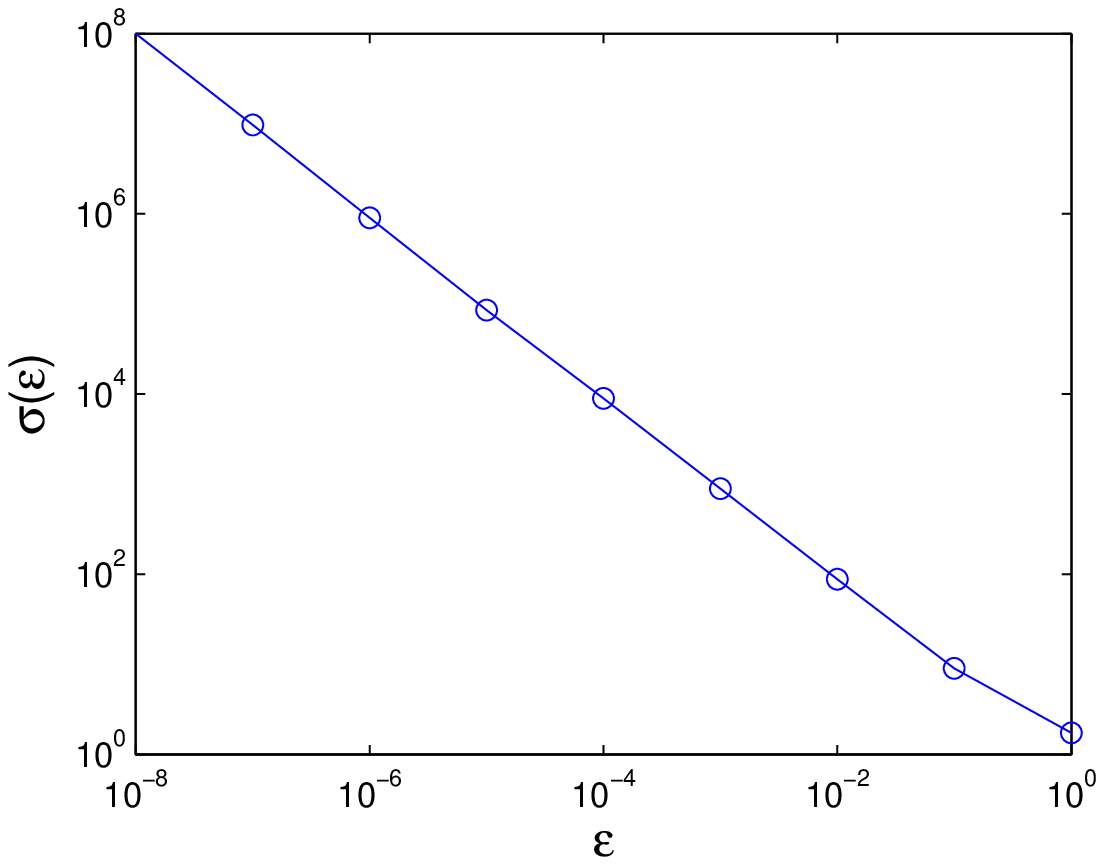}}
\caption{\protect\label{figbvp2} Estimated stiffness ratio of problem
(\ref{prob19}).}
\end{figure}

\begin{exa} Let us consider the following nonlinear problem:
\begin{equation}\label{prob19}
\varepsilon y'' + \exp(y) y' - \frac{\pi}{2} \sin\left( \frac{\pi
t}{2}\right) \exp(2 y) = 0, \ \ \ y(0)=0,\quad y(1)=0.
\end{equation} This problem has a boundary layer at $t=0$ (see
Figure \ref{figbvp2f}). In Figure \ref{figbvp2} we plot an
estimate of the stiffness ratio obtained by considering two
different perturbations of the boundary conditions of the form
$(1, 0)^T$ and $(0,1)^T$. The parameter $\varepsilon$ varies from
$1$ to  $10 ^{-8}$. We see that the (estimated) stiffness
parameter grows like $\varepsilon^{-1}$, as $\varepsilon$ tends to
0.

\end{exa}

\begin{figure}[ht]

\centerline{\includegraphics[width=12cm,height=8cm]{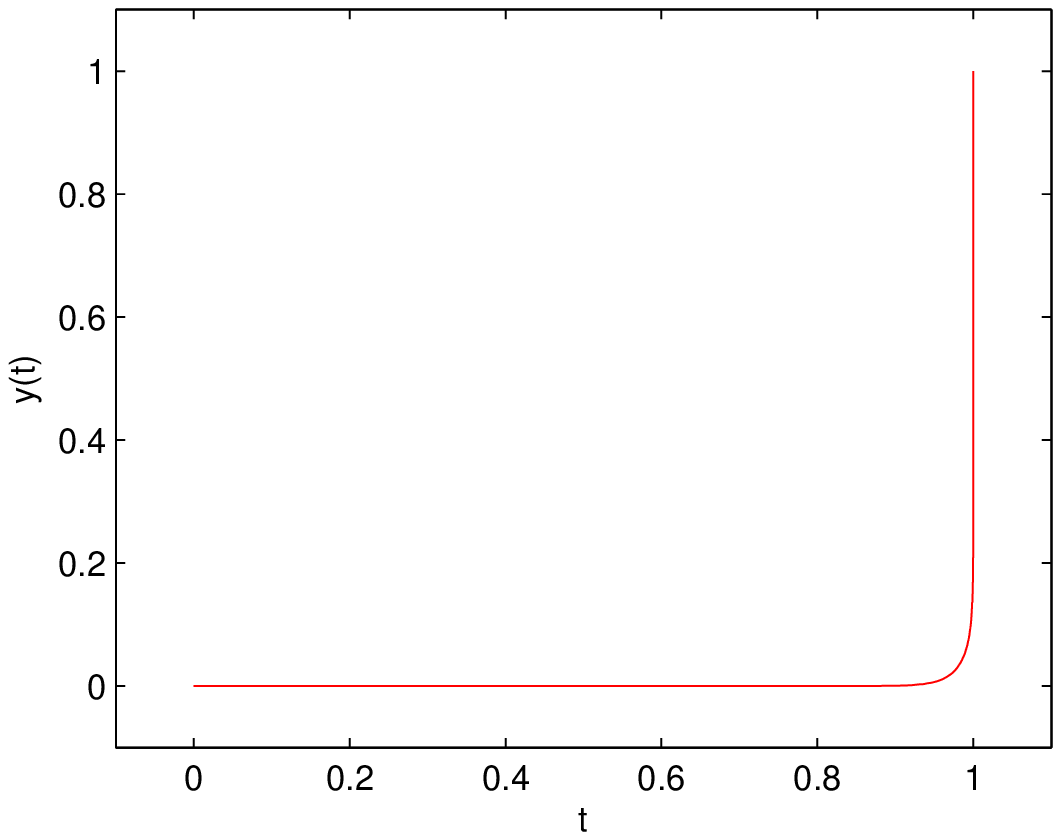}}
\caption{\protect\label{figbvp3f}  Troesch's problem
(\ref{prob6}), $\mu=50$.}
\end{figure}

\begin{figure}[t]

\centerline{\includegraphics[width=12cm,height=8cm]{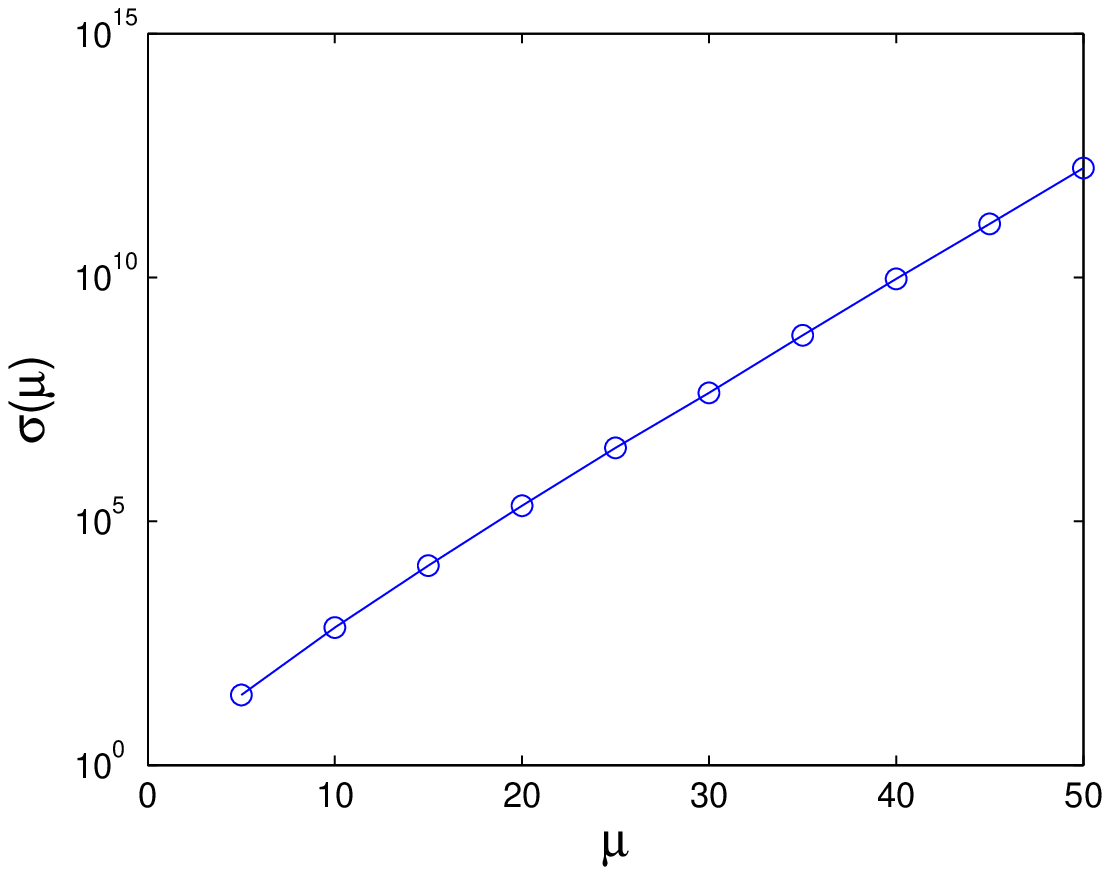}}
\caption{\protect\label{figbvp3} Estimated stiffness ratio of Troesch's
problem (\ref{troesch}).}
\end{figure}

\begin{exa}
Let us consider the nonlinear Troesch problem:
\begin{equation}\label{troesch}
y'' = \mu\, \sinh( \mu y) , \ \ \ y(0)=0,\quad y(1)=1.
\end{equation}
This problem has a boundary layer near $t=1$ (see Figure
\ref{figbvp3f}). In Figure \ref{figbvp3} we plot the estimate of
the stiffness ratio obtained by considering two different
perturbations of the boundary conditions of the form $(1, 0)^T$
and $(0,1)^T$. The parameter $\mu$ is increased from 1 to 50 and,
as expected, the stiffness ratio increases as well: for $\mu=50$,
it reaches the value $1.74 \cdot 10^{12}$.

\end{exa}

\subsection*{Acknowledgements} The authors wish to thank the
reviewers, for their comments and suggestions.

\end{document}